\documentclass[11pt,leqeq]{article}
\usepackage{amssymb,amsthm}
\usepackage{amsmath}
\usepackage{amscd}
\def\Par{{\rm{Par}}}
\def\meet{{\rm{meet}}}
\def\joint{{\rm{joint}}}

\newtheorem{thm}{Theorem}
\newtheorem{lema}[thm]{Lemma}

\newtheorem{defi}[thm]{Definition}

\newcommand{\des}{\displaystyle}
\date{}
\begin{document}
\title{P\'olya Theory for Orbiquotient Sets}
\author{\  Rafael D\'\i az and  H\'{e}ctor Blandin}
\maketitle
\begin{abstract}
Replacing the usual notion of quotient sets by the notion of
orbiquotient sets we obtain a  generalization of P\'olya theory.
The key ingredient of our extended theory is the definition of the
orbicycle  index polynomial which we compute in several examples.
We apply our theory to the study of orbicycles
on orbiquotient sets.\\

\end{abstract}
\section{Introduction}

Assume that a finite group $G$ acts on a finite set $X.$ The
quotient $X/G$ of the action of $G$ on $X$ is a rich and subtle
concept, traditionally $X/G=\{\overline{x}\mid x\in X\}$ where
$\overline{x}=\{gx\mid x\in X\}.$ In recent years it has proven
convenient to modify this notion  in various contexts. For example
one may think of $X/G$ as the groupoid whose set of objects is $X$
and with morphisms given by $X/G(a,b)=\{g\in G\mid ga=b\}$  for
$a, b \in G.$ Following Connes \cite{AConnes} the groupoid $X/G$
is studied with the methods of non-commutative geometry, i.e.
looking at the convolution (incidence) algebra of $X/G.$ Another
approach to quotient sets became rather popular after Vafa and
Witten introduced in \cite{VafaWitten} the so called stringy Euler
numbers. In a nutshell they considered the Euler numbers of
orbiquotient sets
$$X/^{orb}G=\bigsqcup_{\overline{g}\in C(G)}{X^{g}}/{Z(g)},$$
where $C(G)$ is the set of conjugacy classes of $G$, $X^{g}
\subseteq X$ is the set of points fixed by $g\in G$, and $Z(g)$ is
the centralizer of $g$ in $G.$ We will assume that a
representative $g \in G$ has been chosen for each conjugacy class
$\overline{g}$ of $G$. Orbiquotient sets first appeared, see
\cite{as}, in the context of equivariant K-theory in the works of
Atiyah and Segal. The goal of this paper is to bring the notion of
orbiquotient sets into combinatorial waters. Let us provide a
combinatorial motivation for the study of orbiquotient sets
inspired by an analogue topological construction given by
Hirzebruch and H$\ddot{\mbox{o}}$fer in \cite{h}. Consider the set
of $n$-cycles in $X/G,$ i.e., the set of maps
$f:\mathbb{Z}\longrightarrow X/G$ such that $f(k)=f(k+n)$ for
$k\in \mathbb{Z}.$ Suppose we want to lift $f$ to a map
$l:\mathbb{Z}\longrightarrow X$ such that $\pi
\circ l=f$, where $\pi: X \longrightarrow X/G$ is the canonical
projection. The lift $l$ will not be unique, indeed if $l$ is a
lift then $gl$ is another a lift; also we have that $l(k)=gl(k+n)$
for all $k
\in \mathbb{Z}$ and some $g\in G$. Thus the set of lifts of $n$-periodic maps
$\mathbb{Z}\longrightarrow X/G$ may be identified with
$$\{l:\mathbb{Z}\longrightarrow X\mid l(k)=gl(k+n) \mbox{ for } k \in \mathbb{Z}
\mbox{ and some } g \in G\}/G.$$ Inside the later set sits $I(G,X)/G$ the set of constant maps,
where
$$I(G,X)=\{(g,x)\in G \times X \ | \ gx=x\}$$ is the so called inertial set \cite{k}.
The group $G$ acts on $I(G,X)$ as $k(g,x)=(kgk^{-1},kx),$ and it
is not hard to see that
$$I(G,X)/G =\bigsqcup_{\overline{g}\in C(G)}{X^{g}}/{Z(g)}=X/^{orb}G.$$

In this paper we develop orbianalogues for two main results in
elementary combinatorics, the orbit counting lemma and the
Polya-Redfield theorem, see \cite{Bergeron, GCRota}. We fix a
commutative ring $\mathbb{A}$ and consider the category of
$\mathbb{A}$-weighted sets whose objects are pairs $(X,f)$ where
$X$ is a finite set and $f:X\rightarrow\mathbb{A}$ is an arbitrary
map called the weight of $X$. Morphisms between
$\mathbb{A}$-weighted sets are weight preserving bijections. The
cardinality $|X|_{f}$ of a weighted set $(X,f)$ is given by
$$|X|_{f}=\sum_{x\in X}f(x).$$  A finite group $G$ acts on $(X,f)$ if $G$ acts on $X$ and
$f(gx)=f(x)$ for all  $g\in G,$ $x\in X.$ The
Cauchy-Frobenius-Burnside orbit counting lemma gives us a way to
compute $|X/G|_{f}$ as follows:
$$|X/G|_{f}=\frac{1}{|G|}\sum_{g\in G}|X^{g}|_{f}.$$
Suppose now that $G$ is a group of permutations $G\subset S_{m}$.
The cardinality of $X^{m}/G$ is determined by the P\'olya-Redfield
theorem:
$$|X^{m}/G|_{f}=P_{G}(|X|_{f},|X|_{f^{2}},\dots,|X|_{f^{m}}),$$
where $P_{G}$ is the cycle index polynomial of $G$ given by
$$P_{G}(x_{1},x_{2},\dots,x_{m})=\frac{1}{|G|} \sum_{g\in
G}x_{1}^{c_{1}(g)}\dots x_{m}^{c_{m}(g)},$$ and $c_{i}(g)$ is the
number of $g$-cycles of length $i.$ If $X=[n]$ and $f(i)=x_{i}$
for $i\in X$, then directly from the definition of quotient sets
we get that
$$|[n]^{m}/G|_{f}=\sum_{(i_{1},\dots,i_{n})\in\mathbb{N}^{n}}
c_{G}\left({i_{1},\dots,i_{n}}\right)x_{1}^{i_{1}}\dots\
x_{n}^{i_{n}},$$ where $c_{G}\left({i_{1},\dots,i_{n}}\right)$
counts the colorations of $[m]$ with $i_{k}$ elements of color
$k\in[n],$ and two colorations are identified if they are linked
by the action of $G$. The P\'olya-Redfield theorem allows us to
compute the coefficients $c_{G}\left({i_{1},\dots,i_{n}}\right)$
in a different way, namely we have that
$$|[n]^{m}/G|_{f}=P_{G}
\left(\sum_{j=1}^{n}x_{j},\sum_{j=1}^{n}x_{j}^{2},\dots,\sum_{j=1}^{n}x_{j}^{m}\right).$$\\

The rest of this work is organized as follows. In Section 2 we
provide an orbi-analogue of the orbit counting lemma. In Section 3
we provide an orbi-analogue of the P\'olya-Redfield theorem in
full generality, we shall see that lattice of partitions plays a
fundamental role in our presentation. In the remaining sections we
explicitly compute the orbicycle index polynomial for various
groups in increasing order of difficulty. In Section 4 we consider
the case of cyclic groups, and apply it to the study of orbicycles
in orbiquotient sets. In Section 5 we consider the full symmetric
group. In a rather dull fashion we may regard combinatorics as
geometry in dimension zero. It is thus rather interesting when one
can show that the zero dimensional combinatorial case determines
the higher dimensional situation. A theorem of this sort is proved
at the end of Section 5 which provides a strong motivation for the
study of orbiquotient sets. In Section 6 we compute the orbicycle
index polynomial for the dihedral groups.

\section{Orbi-analogue of the orbit counting lemma}

If $S\subseteq G$ and $G$ acts on $X$, then we set $X^{S}=\{x\in
X\mid gx=x \ \textrm{for\ } g\in S\}.$ Also we let $\langle
g_{1},\dots,g_{n}\rangle$ be the subgroup of $G$ generated by
$\{g_{1},\dots,g_{n}\}\subset G.$

\begin{defi}\label{DEFI4}{\em
The orbiquotient of $X$ by the action of $G$ is the set given by
$$X/^{orb}G=\bigsqcup_{\overline{g}\in C(G)}{X^{g}}/{Z(g)}.$$}
\end{defi}
The orbiquotient $X/^{orb}G$ is well defined up to canonical
bijections. Indeed if $h=kgk^{-1}$ then the map
$\psi:X^{g}\rightarrow X^{h}$ given by $\psi(x)=kx$ induces a
bijection $$\psi:{X^{g}}/{Z(g)}\rightarrow{X^{h}}/{Z(h)}.$$ If $G$
acts on a weighted set $(X,f)$ then $X/^{orb}G$ is also weighted:
$X^{g}$ is weighted by $f\mid_{X^{g}}$ and ${X^{g}}/{Z(g)}$ is
weighted by $f(\overline{x})=f(x)$ for
$\overline{x}\in{X^{g}}/{Z(g)}.$ Our next result is the
orbi-analogue of the orbit counting lemma, let us first introduce
a notation that will be used repeatedly
$$P(G)=\{(\overline{g},h)\ | \ \overline{g} \in C(G) \mbox{ and }
h \in Z(g\}.$$

\begin{thm}\label{T4}{\em If $G$ acts on $(X,f),$ then the cardinality of $X/^{orb}G$ is given by
$$\left|X/^{orb} G\right|_{f}=\displaystyle{\frac{1}{|G|}\sum_{(\overline{g},h)
\in P(G)}|\overline{g}|\left|X^{\langle
g,h\rangle}\right|_{f}}.$$}
\end{thm}
The proof of this result is quite simple:

\begin{eqnarray*} \left|X/^{orb} G\right|_{f}
&=&\sum_{\overline{g}\in C(G)}\left|{X^{g}}/{Z(g)}\right|_{f}\\
&=&\sum_{\overline{g}\in C(G)}\frac{|\overline{g}|}{|G|}\sum_{h\in Z(g)}|X^{g}\cap X^{h}|_{f}\\
&=& \frac{1}{|G|}\sum_{(\overline{g},h)\in
P(G)}|\overline{g}|\left|X^{\langle g,h\rangle}\right|_{f}.
\end{eqnarray*}

\section{Orbi-analogue of P\'olya-Redfield theorem}

Let $\Par(X)$ be lattice of partitions of $X.$  The minimal and
maximal elements of $\Par(X)$ are $\{\{x\}\mid x\in X\}$ and
$\{X\},$ respectively.  The $\joint$ $\pi\vee\rho$ of partitions
$\pi$ and $\rho$ is defined by demanding that $i,j\in X$ belong to
a block of $\pi\vee\rho$ if there exists a sequence
$i=a_{0},a_{1},\dots,a_{n}=j,$ such that for  $0\leq i
\leq n-1$ either $a_{i}$ and $a_{i+1}$ belong to a block in
$\pi,$ or $a_{i}$ and $a_{i+1}$ belong a block in $\rho.$ The
$\meet$ of partitions $\pi$ and $\rho$ is $\pi\wedge\rho=\{B\cap
C\mid B\in\pi, C\in\rho\rm,\it \ B\cap{C}\neq \mathrm{0} \}.$ Let
the group $G$ act on a set $X$ with $n$-elements. Each $g \in G$
induces a partition $C(g)$ on $X$ such that
$C(g)=\bigsqcup_{i=1}^{n}C_{i}(g),$  where
$C_{i}(g)=\{g\hbox{-}\rm cycles\ on\ X\ of\ length\ \it i\}$ for
$1\leq i\leq n.$ We use the notation $c(g)=|C(g)|$ and
$c_{i}(g)=|C_{i}(g)|$ for $1\leq i\leq n.$  If $\pi$ is a
partition of $X$ we let $b_k(\pi)$ be the number of blocks of
$\pi$ of cardinality $k.$

\begin{defi}{\em
The orbicycle index polynomial $P^{orb}_{G}(x_{1},x_{2},\dots)\in
\mathbb{Q}[x_{1},x_{2},\dots]$ is given by
$$P_{G}^{\
orb}(x_{1},x_{2},\dots)=\frac{1}{|G|}\sum_{(\overline{g},h)\in
P(G)}|\overline{g}|x^{C(g)\vee C(h)},$$  where
$$x^{C(g)\vee C(h)}=\prod_{k\geq 1} x_{k}^{b_{k}(C(g)\vee C(h))}.$$}
\end{defi}

If $G \subseteq S_m$ then $G$ acts on $X^{m}$. Suppose that
$g,h\in G$ commute, then $i,j\in[m]$ belong to the same block of
$C(g)\vee C(h)$ if and only if there exist $a,b\in\mathbb{Z}$ such
that $j=(g^{a}h^{b})(i).$ It is easy to check  that  $f \in X^{m}$
is fixed by $g$ and $h$ if and only if $f$ is constant on each
block of $C(g)\vee C(h).$

\begin{thm}{\em
Let $(X,f)$ be an $\mathbb{A}$-weighted set and $G \subseteq S_m.$
The cardinality of $X^{m}/^{orb} G$ is given by
\[|X^{m}/^{orb} G|_{f}=P_{G}^{\
orb}(|X|_{f},|X|_{f^{2}},\dots).\]}
\end{thm}

\begin{proof}

\begin{eqnarray*}
|X^{m}/^{orb}
G|_{f}&=&\des{\frac{1}{|G|}\sum_{(\overline{g},h)\in P(G)}|\overline{g}|\left|
\left(X^{[m]}\right)^{g}\cap\left(X^{[m]}\right)^{h}\right|_{\tilde{f}}}\\
&=&\des{\frac{1}{|G|}\sum_{(\overline{g},h)\in
P(G)}|\overline{g}|\sum_{\alpha:[m]\rightarrow X,
\alpha\circ g=\alpha\circ h=\alpha }\prod_{x\in[m]}f(\alpha(x))}\\
&=&\des{\frac{1}{|G|}\sum_{(\overline{g},h)\in
P(G)}|\overline{g}|\sum_{\alpha:C(g)\vee C(h)\rightarrow X}
\prod_{B\in C(g)\vee C(h)}\prod_{x\in B}f(\alpha(x))}\\
&=&\des{\frac{1}{|G|}\sum_{(\overline{g},h)\in P(G)}|\overline{g}|\prod_{B\in C(g)\vee C(h)}\sum_{y\in X}f(y)^{|B|}}\\
&=&\des{\frac{1}{|G|}\sum_{(\overline{g},h)\in P(G)}|\overline{g}|\prod_{k\geq 1}|X|_{f^{k}}^{b_{k}(C(g)\vee C(h))}}\\
&=&P_{G}^{\ orb}(|X|_{f},|X|_{f^{2}},\dots).
\end{eqnarray*}
\end{proof}

Let $X=[n]$ and $f(i)=x_i$, then one can check directly from the
definition that
$$|[n]^{m}/^{orb} G|_{f}=\sum_{(i_{1},\dots,i_{n}) \in \mathbb{N}^n}
c^{orb}_{G}\left({i_{1},\dots,i_{n}}\right)x_{1}^{i_{1}}\dots\
x_{n}^{i_{n}},$$ where
$c^{orb}_{G}\left({i_{1},\dots,i_{n}}\right)$ counts colorations
$c$ of $[m]$ with colors in $[n]$ such that:
\begin{itemize}
\item  There are $i_{k}$ elements in $[m]$ of color $k\in[n].$
\item $c$ is $g$-invariant for some $\overline{g} \in C(G).$
\item Two $g$-invariant colorations are identified if
they can be linked by the action of $Z(g).$
\end{itemize}

The orbi-analgogue of the P\'olya-Redfield gives us another way to
compute the coefficients
$c^{orb}_{G}\left({i_{1},\dots,i_{n}}\right)$, namely we have that
$$|[n]^{m}/^{orb} G|_{f}=P^{\ orb }_{G}
\left(\sum_{j=1}^{n}x_{j},\sum_{j=1}^{n}x_{j}^{2},\dots,\sum_{j=1}^{n}x_{j}^{n}\right).$$

\section{Orbicycle index polynomial of $\mathbb{Z}_{n}$}

Let $\mathbb{N_{+}}$ be the set of positive integers and let
$(x_{1},x_{2},\dots,x_{k})$ be the greatest common divisor of
$x_{1},x_{2},\dots,x_{k}\in\mathbb{N}^{+}$. The cyclic group with
$n$-elements is denoted by $\mathbb{Z}_{n}=\{1,2,\dots,n\}$. For $n,k \in
\mathbb{N_{+}}$ we define an equivalence
relation on $\mathbb{Z}_{n}^{k}$ as follows: $x$ and $y$ are
equivalent if and only if $(x,n)=(y,n).$  It is easy to verify
that
$\mathbb{Z}_{n}^{k}=\bigsqcup_{d|n}\{x\in\mathbb{Z}_{n}^{k}\mid
(x,n)=d\},$ and thus we have
$$n^{k}=\sum_{d\mid n}|\{x\in\mathbb{Z}_{n}^{k}\mid (x,n)=d\}|.$$
For $n\in \mathbb{N_{+}}$ the Jordan totient function $J_{k}$, see
\cite{TomMikeApostol}, is given by $J_{k}(n)=\left|\{x\in\mathbb{Z}_{n}^{k}\mid
(x,n)=1\}\right|.$  For each  $d|n$ we have
$J_{k}\left(\frac{n}{d}\right)=\left|\{x\in\mathbb{Z}_{n}^{k}\mid
(x,n)=d\}\right|,$ therefore we get that $n^{k}=\sum_{d\mid
n}J_{k}(d).$ By the M$\ddot{\hbox{o}}\rm bius$ inversion formula
$J_{k}(n)=\sum_{d\mid n}\mu\left(\frac{n}{d}\right)d^{k},$ thus
$J_{k}(p^{r})=p^{kr}-p^{k(r-1)},$ for $p$ prime, and for arbitrary
integer $n=p_1^{\alpha_1} \dots p_r^{\alpha_r}$ where $p_{1},\dots,p_{r}$ are distinct prime numbers, we get
$$J_{k}(n)=n^{k}\left(1-\frac{1}{p_{1}^{k}}\right)\dots\left(1-\frac{1}{p_{r}^{k}}\right).$$
We shall need the following property, an easy consequence of the
previous considerations, of the Jordan totient function:
$$\sum_{x\in\mathbb{Z}_{n}^{k}}f((x,n))
=\sum_{d\mid n}J_{k}\left(\frac{n}{d}\right)f(d),$$  for any $f:\{d: d\mid
n\}\longrightarrow\mathbb{A}.$ Recall that if
$x,y\in\mathbb{Z}_{n}\subseteq S_n,$ then $|C(x)\vee
C(y)|=(x,y,n),$ and all blocks in $C(x)\vee C(y)$ are of
cardinality $\frac{n}{(n,x,y)}.$ Indeed if $x\in\mathbb{Z}_{n},$
then $\frac{\mathbb{Z}_{n}}{(x)}\cong
\mathbb{Z}_{(n,x)},$ since for $a,b\in\mathbb{Z}_{n}$ we have that
$a=b$ in $\mathbb{Z}_{n}/(x)$ if and only if $a=b$ mod $(n,x).$
Thus there are $(n,x)$ blocks in $\mathbb{Z}_{n}/(x)$ all of them
of cardinality $\frac{n}{(n,x)}.$ Similarly if
$x,y\in\mathbb{Z}_{n},$ then $a,b\in\mathbb{Z}_{n}$ are in the
same block of $C(x)\vee C(y)$ if and only if there exist
$r,s\in\mathbb{Z}$ such that $b=a+rx+sy$ mod $n$, or equivalently
$a=b$ mod $(n,x,y).$

\begin{thm}{\em
$$P_{\mathbb{Z}_{n}}^{\
orb}(y_{1},y_{2},\dots,y_{n})=\frac{1}{n}\sum_{d\mid
n}J_{2}(d)y_{d}^{\frac{n}{d}}.$$}
\end{thm}

\begin{proof}
\begin{eqnarray*}
P_{\mathbb{Z}_{n}}^{\
orb}(y_{1},y_{2},\dots,y_{n})&=&\frac{1}{n}\sum_{(\overline{g},h)\in
P(\mathbb{Z}_{n})}|\overline{g}|\prod_{k\geq
1}y_{k}^{b_{k}(C(g)\vee C(h))}\\
&=&\frac{1}{n}\sum_{(x,y)\in\mathbb{Z}_{n}\times\mathbb{Z}_{n}}\prod_{k\geq
1}y_{k}^{b_{k}(C(x)\vee C(y))}\\
&=&\frac{1}{n}\sum_{(x,y)\in\mathbb{Z}_{n}\times\mathbb{Z}_{n}}y_{\frac{n}{(x,y,n)}}^{(x,y,n)}\\
&=&\frac{1}{n}\sum_{d\mid
n}J_{2}\left(d\right)y_{d}^{\frac{n}{d}}.
\end{eqnarray*}
\end{proof}

The coefficients
$c_{\mathbb{Z}_n}\left({i_{1},\dots,i_{m}}\right)$ are computed in
\cite{LeonidZimmels}. As a corollary of the previous result we obtain that
$$c_{\mathbb{Z}_n}^{orb}\left({i_{1},\dots,i_{m}}\right)=\frac{1}{n}\sum_{d\mid
(i_{1},\dots,i_{m})}J_{2}(d)\frac{\left(\frac{i_{1}}{d}+\dots+\frac{i_{m}}{d}\right)!}
{\frac{i_{1}}{d}!\dots\frac{i_{m}}{d}!}.$$

Let $p$ be a prime number and $r\in\mathbb{N}^{+},$ necklaces
without a clasp with $p$ beads and $r$ colors may be identified
with the set $C_{p}([r])=[r]^{p}/\mathbb{Z}_{p}.$ As explained in
\cite{DZeilberger}  we have that
$$|C_{p}([r])|=P_{\mathbb{Z}_{p}}(r,r,\dots)=\frac{1}{p}\sum_{d\mid
p}\varphi(d)r^{\frac{p}{d}}=r+\frac{r^{p}-r}{p}.$$

\begin{defi}{\em
The set $C_{n}^{orb}(X)$ of orbi $n$-cycles in $X$ is given by
 $C^{orb}_{n}(X)=X^{n}/^{orb}\mathbb{Z}_{n}.$}
\end{defi}

In analogy with the example above we define the set of
orbi-necklaces without a clasp with $p$ beads and $r$ colors to be
$C^{orb}_{p}([r])=[r]^{p}/^{orb}\mathbb{Z}_{p}.$ Its cardinality
is given by
$$|C^{orb}_{p}([r])|=P_{\mathbb{Z}_{p}}^{\ orb}(r,r,\dots)=\frac{1}{p}\sum_{d\mid
p}J_{2}(d)r^{\frac{p}{d}}=rp+\frac{r^{p}-r}{p}.$$ Next couple of
results count explicitly the number of orbicycles in orbiquotient
sets.

\begin{thm}{\em
If $G$ acts on $(X,f)$ then
$$|C^{orb}_{n}(X/^{orb}G)|_{f}=\frac{1}{n}\sum_{
\alpha:[\frac{n}{d}]\rightarrow P(G)}\frac{J_{2}(d)}{|G|^{\frac{n}{d}}}\prod_{i=1}^{\frac{n}{d}}
|\pi_{C(G)}(\alpha(i))||X^{\alpha(i)}|_{f^{d}}.$$}
\end{thm}
\begin{proof}
\begin{eqnarray*}
|C^{orb}_{n}(X/^{orb} G)|_{f}&=& P_{\mathbb{Z}_{n}}^{\ orb}(|X/^{orb} G|_{f},|X/^{orb} G|_{f^{2}},\dots)\\
&=&\frac{1}{n}\sum_{d\mid n}J_{2}(d)\left(|X/^{orb} G|_{f^{d}}\right)^{\frac{n}{d}}\\
&=&\frac{1}{n}\sum_{d\mid
n}J_{2}(d)\left(\frac{1}{|G|}\des{\sum_{(\overline{g},h)\in
P(G)}|\overline{g}||X^{\langle g,h \rangle}|_{f}}\right)^{\frac{n}{d}}\\
&=& \frac{1}{n}\sum_{d\mid
n}\frac{J_{2}(d)}{|G|^{\frac{n}{d}}}\sum_{\alpha:[\frac{n}{d}]\rightarrow
P(G)}\prod_{i=1}^{\frac{n}{d}}|\pi_{C(G)}(\alpha(i))|
|X^{\langle \alpha(i)\rangle}|_{f^{d}}\\
&=&\frac{1}{n}\sum_{\alpha:[\frac{n}{d}]\rightarrow
P(G)}\frac{J_{2}(d)}{|G|^{\frac{n}{d}}}\prod_{i=1}^{\frac{n}{d}}
|\pi_{C(G)}(\alpha(i))||X^{\alpha(i)}|_{f^{d}}.
\end{eqnarray*}
\end{proof}

\begin{thm}{\em
Let $(X,f)$ be be an $\mathbb{A}$-weithed set and  $G \subseteq
S_m$, then we have
$$|C^{orb}_{n}(X^{[m]}/^{orb}\ G)|_{f}=\frac{1}{n}\sum_{ \alpha:[\frac{n}{d}]\rightarrow P(G)}
\frac{J_{2}(d)}{|G|^{\frac{n}{d}}}\des{\prod_{i=1}^{\frac{n}{d}}|X|_{f^{d}}^{b(\alpha(i))}}$$
where $|X|_{f^{d}}^{b(\alpha(i))}=\prod_{k\geq
1}|X|_{f^{dk}}^{b_{k}(C(\pi_{1}(\alpha(i)))\vee
C(\pi_{2}(\alpha(i))))}$}
\end{thm}
\begin{proof}
\begin{eqnarray*}
|C^{orb}_{n}(X^{[m]}/^{orb} G)|_{f} &=& P_{\mathbb{Z}_{n}}^{\
orb}(|X^{[m]}/^{orb} G|_{f},|X^{[m]}/^{orb} G|_{f^{2}},\dots)\\
&=& \frac{1}{n}\sum_{d\mid n}J_{2}(d)\des{P_{G}^{\ orb}(|X|_{f^{d}},|X|_{f^{2d}},\dots)}^{\frac{n}{d}}\\
&=& \frac{1}{n}\sum_{d\mid
n}J_{2}(d)\left(\des{\frac{1}{|G|}\sum_{(\overline{g},h)\in
P(G)}\prod_{k\geq 1}
|X|_{f^{dk}}^{b_{k}(C(g)\vee C(h))}}\right)^{\frac{n}{d}}\\
&=& \frac{1}{n}\sum_{d\mid
n}\frac{J_{2}(d)}{|G|^{\frac{n}{d}}}\des{\sum_{\alpha:[\frac{n}{d}]
\rightarrow P(G)}\prod_{i=1, k\geq 1}|X|_{f^{dk}}^{b_{k}(C(\pi_{1}(\alpha(i)))\vee C(\pi_{2}(\alpha(i))))}}\\
&=&
\frac{1}{n}\sum_{\alpha:[\frac{n}{d}]\rightarrow
P(G)}\frac{J_{2}(d)}{|G|^{\frac{n}{d}}}\des{\prod_{i=1}|X|_{f^{d}}^{b(\alpha(i))}}.
\end{eqnarray*}
\end{proof}
Using similar methods one can count cycles on orbiquotient sets:
$$|C_{n}(X/G)|_{f}=\frac{1}{n} \sum_{\alpha:[\frac{n}{d}]\rightarrow G}\frac{J_{2}(d)}{|G|^{\frac{n}{d}}}
\prod_{i=1}^{\frac{n}{d}}|X^{\alpha(i)}|_{f^{d}},$$
and
$$|C_{n}(X^{m}/G)|_{f}=\frac{1}{n}\sum_{\alpha:[\frac{n}{d}]\rightarrow G}
\frac{J_{2}(d)}{|G|^{\frac{n}{d}}}
\des{\prod_{i=1}^{\frac{n}{d}}|X|_{f^{d}}^{c(\alpha(i))}},$$
where $|X|_{f^{d}}^{c(g)}=\des{\prod_{k\geq
1}|X|_{f^{dk}}^{c_{k}(g)}}$ for $g\in G.$

\section{Orbicyle index polynomial of $S_{n}$}

A partition of depth $k$, denoted by $\alpha\vdash_{k}n$,  of
$n\in\mathbb{N}_{+}$ is a map
$\alpha:\left(\mathbb{N}_{+}\right)^{k}\rightarrow\mathbb{N}$ such
that $$\sum_{(i_{1},\dots,i_{k})\in \mathbb{N}_{+}^k}i_{1}\dots
i_{k}\alpha(i_{1},\dots,i_{k})=n.$$  A partition of depth $1$ is a
partition in the usual sense. To each partition $\alpha$ we
associate a canonical permutation of $[n]$ whose cycle structure
is determined by $\alpha.$ Keeping this correspondence in mind one
can check that if $\alpha\vdash n$ then
\begin{enumerate}
\item $Z(\alpha)$ is isomorphic to $\prod_{i=1}^{n}\mathbb{Z}_{i}^{\alpha_{i}}\rtimes
S_{\alpha_{i}}.$
\item If $h\in Z(\alpha),$ then $b_{k}(C(\alpha)\vee C(h))=\sum_{d\mid
k}c_{\frac{k}{d}}(\pi_{d}(h))$ where $\pi_{d}(h)$ is the
projection of $h$ into $\prod_{i=1}^{n}S_{\alpha_{i}}.$
\end{enumerate}

\begin{thm}{\em
$$P^{\ orb}_{S_{n}}(y_{1},y_{2},\dots,y_{n})=
\sum_{\beta{\vdash_{2}}n}\prod_{(i,j,k)\in[n]
\times[\alpha_{i}]\times[n]}\frac{y_{k}^{\sum_{d\mid
k}\beta\left(d,\frac{k}{d}\right)}}{j^{\beta(i,j)}\beta(i,j)!}.$$}
\end{thm}
\begin{proof}
By the previous remarks we have
\begin{eqnarray*}
P^{orb}_{S_{n}}(y_{1},y_{2},\dots,y_{n})&=&\frac{1}{n!}\sum_{(\overline{g},h)\in
P(S_{n})}|\overline{g}|\prod_{k=1}^{n}y_{k}^{b_{k}(C(g)\vee
C(h))}\\
&=&\frac{1}{n!}\sum_{\alpha\vdash
n}\frac{n!}{\prod_{i=1}^{n}i^{\alpha_{i}}\alpha_{i}!}
\sum_{h\in
Z(\alpha)}\prod_{k=1}^{n}y_{k}^{{b}_{k}(C(\alpha)\vee
C(h))}\\
&=&\sum_{\alpha\vdash n}\frac{1}{\prod_{i=1}^{n}\alpha_{i}!}
\sum_{h\in\prod_{i=1}^{n}S_{\alpha_{i}}}
\prod_{k=1}^{n}y_{k}^{\sum_{d\mid
k}c_{\frac{k}{d}}(\pi_{d}(h))}\\
&=&\sum_{\alpha\vdash
n}\frac{1}{\prod_{i=1}^{n}\alpha_{i}!}\sum_{{\beta_{i}
\vdash\alpha_{i}}}\frac{\prod_{i=1}^{n}\alpha_{i}!}{\prod_{(i,j)\in
[n]\times[\alpha_{i}]}j^{\beta(i,j)}\beta(i,j)!}\prod_{k=1}^{n}y_{k}^{
\sum_{d\mid k}\beta_{\frac{k}{d}}(d)}\\
&=&\sum_{\beta{\vdash_{2}}n}\prod_{(i,j,k)\in[n]
\times[\alpha_{i}]\times[n]}\frac{y_{k}^{\sum_{d\mid
k}\beta\left(d,\frac{k}{d}\right)}}{j^{\beta(i,j)}\beta(i,j)!}.
\end{eqnarray*}
Above we used the fact that $\sum_{d\mid
k}c_{\frac{k}{d}}(\pi_{d}(h))$ depends only on the cycle structure
of $\pi_{d}(h).$
\end{proof}

The orbicycle index polynomial can be use to compute the even
dimensions of the orbifold cohomology groups for global orbifolds
of the form $M^{n}/^{\textit{orb}}G,$ where $M$ is a compact
smooth manifold, and $G\subset S_n.$   For simplicity we only
consider cohomology in even dimensions. The orbifold cohomology is
defined as follows:
$$\mathrm{H}^{orb}\left(M^{n}/G\right)=\bigoplus_{\bar{g}\in
C(G)}\mathrm{H}\left(\left(M^{n}\right)^{g}\right)^{Z(g)}.$$

The following result is a direct consequence of the
characterization of the centralizer of permutations previously
discussed.
\begin{lema}
$$\mathrm{H^{\textit{orb}}}\left(M^{n}/G\right)=
\bigoplus_{\bar{g}\in C(G)}
\bigotimes_{i}\left({\mathrm{H}\left(M\right)}^{\otimes
c_{i}(g) }\right)^{S_{c_{i}(g)}}.$$
\end{lema}
\begin{proof}
\begin{align*}
\mathrm{H^{orb}}\left(M^{n}/G\right)=&\bigoplus_{\bar{g}\in
C(G)}\mathrm{H}\left(\left(M^{n}\right)^{g}\right)^{Z(g)}\\
=&\bigoplus_{\bar{g}\in
C(G)}\mathrm{H}(\prod_{i}M^{c_{i}(g)})^{Z(g)}\\
=&\bigoplus_{\bar{g}\in C(G)}
\bigotimes_{i}\left({\mathrm{H}\left(M\right)}^{\otimes
c_{i}(g) }\right)^{S_{c_{i}(g)}}.
\end{align*}
\end{proof}

Assume that we are given a finite basis $X$ for $\mathrm{H}(M)$,
then we have the following result.

\begin{thm}{\em $$\dim\left(\mathrm{H}^{orb}\left(M^{n}/S_{n}\right)\right)=
 P_{G}^{\ orb}(|X|,\dots,|X|).$$}
\end{thm}
\begin{proof}
\begin{align*}
\dim\left(\mathrm{H}^{orb}\left(M^{n}/S_{n}\right)\right)
&=\sum_{\bar{g}\in C(G)}\prod_{i}\left|X^{c_{i}(g)}/S_{c_{i}(g)}\right|\\
&=\left|X^{n}/^{orb}S_{n}\right|\\
&=P_{G}^{\ orb}(|X|,\dots,|X|).
\end{align*}
\end{proof}

Notice that above we use the trivial weight on $X$; using a
generic weight we obtain further information on the orbifold
cohomology groups. Theorem $12$ gives a combinatorial
interpretation for the orbifold cohomology groups, however we do
not have a combinatorial interpretation for the orbifold product
introduced by Chen and Ruan in \cite{c}. Until recently this
problem seemed hopeless, however the alternative description of
the Chen-Ruan product introduced by Jarvis, Kauffman and Kimura in
\cite{k} could pave the way for such a combinatorial
understanding. Theorem $11$ suggests the possibility of
constructing, along the lines of \cite{QSF}, an orbi-analogue for
the symmetric functions. This issue deserves further research.

\section{Orbicycle index polynomial of $D_{n}$}

The generators $\rho$ and $\tau$ of the dihedral group
$D_{n}=\{e,\rho,\dots,\rho^{n-1},\tau,\dots,\tau\rho^{n-1}\}$ are
such that $\rho^{n}=e$, $\tau^{2}=e$ and
$\tau\rho=\rho^{n-1}\tau.$ The conjugacy classes of the dihedral
groups are described in the following tables, see
\cite{DancerIsaacLinks}. For $n$ odd there are $\frac{n+3}{2}$
conjugacy classes organized in three families
\begin{center}
\begin{tabular}{|c|c|c|}
  \hline
  Conjugacy class & Representative & Centralizer subgroup \\
  \hline
  $\{e\}$ & $e$ & $D_{n}$ \\
  \hline
  $\{\rho^{i},\rho^{-i}\}$\ for\ $1\leq i\leq\frac{n-1}{2}$ & $\rho^{i}$ & $\{\rho^{i}\mid 0\leq i<n\}$ \\
  \hline
  $\{\rho^{i}\tau\mid 0\leq i<n\}$ & $\tau$  & $\{e,\tau\}$ \\
  \hline
\end{tabular}
\end{center}
For $n$ even there are $\frac{n+6}{2}$ conjugacy classes organized
in five families
\begin{center}
\begin{tabular}{|c|c|c|}
  \hline
  Conjugacy class & Representative & Centralizer subgroup \\
  \hline
  $\{e\}$ & $e$ & $D_{n}$ \\
  \hline
  $\{\rho^{\frac{n}{2}}\}$ & $\rho^{\frac{n}{2}}$ & $D_{n}$ \\
  \hline
  $\{\rho^{i},\rho^{-i}\},$ \ $1\leq i<\frac{n}{2}$ & $\rho^{i}$ & $\{\rho^{i}\mid 0\leq i<n\}$ \\
  \hline
  $\{\rho^{2i}\tau\mid 0\leq i<\frac{n}{2}\}$ & $\tau$  & $\{e,\tau,\rho^{\frac{n}{2}},\rho^{\frac{n}{2}}\tau\}$ \\
  \hline
  $\{\rho^{2i+1}\tau\mid 0\leq i<\frac{n}{2}\}$ & $\rho\tau$ & $\{e,\rho\tau,\rho^{\frac{n}{2}},\rho^{\frac{n}{2}+1}\tau\}$ \\
  \hline
\end{tabular}
\end{center}
For real numbers $a_{1}\geq 1,\dots,a_{k}\geq 1$, and $n\in\mathbb{N}^{+}$, see \cite{TomMikeApostol}, we set 
$$\varphi(a_{1},\dots,a_{k},n)=|\{x\in\mathbb{Z}_{n}^{k}\mid\ 1\leq x_i\leq a_i,\forall i\in[k],\ (x,n)=1\}|.$$
For $x\in\mathbb{R}$ we denote by $\lfloor x\rfloor$ the floor function of $x$, that is the largest integer no greater than $x$. We use the fact that for a real number $x\geq 1$ and $n\in\mathbb{N}^{+}$ we have 
$$\left\lfloor\frac{x}{n}\right\rfloor=\{a\in\mathbb{Z}:\ 1\leq a\leq x,\ n\mid a\}$$
Proceeding as in the previous sections one can show that

$$\varphi\left(\frac{a_{1}}{d},\dots,\frac{a_{k}}{d},\frac{n}{d}\right)=|\{x\in\mathbb{Z}_{n}^{k}\mid\ 1\leq x_i\leq a_i,\forall i\in[k],\ (x,n)=d\}|.$$

and thus 

\begin{equation}
\sum_{d\vert n}\varphi\left(\frac{a_1}{d},\dots, \frac{a_k}{d},\frac{n}{d}\right)
=\lfloor a_1\rfloor \dots\lfloor a_k\rfloor
\end{equation}

and this implies

$$\varphi(a_{1},\dots,a_{k},n)=\sum_{d\mid n}\mu\left(d\right)
\left\lfloor\frac{a_{1}}{d}\right\rfloor\dots \left\lfloor\frac{a_{k}}{d}\right\rfloor,$$ 

\noindent and also for $f:\{d: d\mid
n\}\rightarrow\mathbb{A}$ an arbitrary map we have that
$$\sum_{x\in\mathbb{Z}_{a_{1}}
\times\dots\times\mathbb{Z}_{a_{k}}}f\left( (x,n) \right)=\sum_{d\mid
n}\varphi\left(\frac{a_1}{d},\dots, \frac{a_k}{d},\frac{n}{d}\right)f(d).$$



\begin{thm}{\em Let $n \in \mathbb{N}_+$ be odd. The orbicycle index polynomial of $D_{n}$  is given by
$$P_{D_{n}}^{\ orb}=\frac{1}{n}\sum_{d\mid
n}\varphi\left(\frac{n-1}{2d},\frac{n-1}{d},\frac{n}{d}\right)x_{d}^{\frac{n}{d}}-\frac{1}{2}P_{\mathbb{Z}_{n}}
+\frac{3}{2}x_{1}x_{2}^{\frac{n-1}{2}}$$ 
}
\end{thm}
\begin{proof}
\begin{eqnarray*}
P_{D_{n}}^{\
orb}(x_{1},\dots,x_{n})&=&\frac{1}{2n}\sum_{(\overline{g},h)\in
P(D_{n})}|\overline{g}|\prod_{k=1}^{n}x_{k}^{b_{k}(C(g)\vee
C(h))}\\
&=&P_{D_{n}}+\frac{1}{n}\sum_{i=1}^{\frac{n-1}{2}}
\sum_{j=0}^{n-1}\prod_{k=1}^{n}x_{k}^{b_{k}(C(\rho^{i})\vee
C(\rho^{j}))}+x_{1}x_{2}^{\frac{n-1}{2}}\\
&=&P_{D_{n}}+\frac{1}{n}\sum_{i=0}^{\frac{n-1}{2}}
\sum_{j=0}^{n-1}\prod_{k=1}^{n}x_{k}^{b_{k}(C(\rho^{i})\vee
C(\rho^{j}))}
-P_{\mathbb{Z}_{n}}+x_{1}x_{2}^{\frac{n-1}{2}}\\
&=&\frac{1}{2}P_{\mathbb{Z}_{n}}+\frac{3}{2}x_{1}x_{2}^{\frac{n-1}{2}}+
\frac{1}{n}\sum_{i=0}^{\frac{n-1}{2}}
\sum_{j=0}^{n-1}\prod_{k=1}^{n}x_{k}^{b_{k}(C(\rho^{i})\vee
C(\rho^{j}))}
-P_{\mathbb{Z}_{n}}\\
&=&\frac{1}{n}\sum_{d\mid
n}\varphi\left(\frac{n-1}{2d},\frac{n-1}{d},\frac{n}{d}\right)x_{d}^{\frac{n}{d}}-\frac{1}{2}P_{\mathbb{Z}_{n}}
+\frac{3}{2}x_{1}x_{2}^{\frac{n-1}{2}}.
\end{eqnarray*}
\end{proof}

Recall that $\tau$ and $\rho$ are given $\tau(x)=3-x$ and
$\rho^{r}(x)=x+r.$ Our next table gives the equivalence class of
$x\in\mathbb{Z}_{n}$ under five different equivalence relations.

\begin{center}
\begin{tabular}{|c|c|}
\hline
Partition & Equivalence class \\
\hline
$C(\tau)\vee C(\rho^{\frac{n}{2}})$ &
$\left\{x,3-x,x+\frac{n}{2},3-x+\frac{n}{2}\right\}$  \\
\hline
$C(\tau)\vee C(\rho^{\frac{n}{2}}\tau)$ & $
\left\{x,3-x,3-x+\frac{n}{2},x+\frac{n}{2}\right\}$  \\
  \hline
$C(\rho\tau)\vee C(\rho^{\frac{n}{2}})$ & $
\left\{x,4-x,x+\frac{n}{2},4-x+\frac{n}{2}\right\}$ \\
  \hline
$C(\rho\tau)\vee C(\rho^{\frac{n}{2}+1}\tau)$ & $
\left\{x,4-x,4-x+\frac{n}{2},x+\frac{n}{2}\right\}$ \\
  \hline
$C(\rho^{\frac{n}{2}})\vee C(\tau\rho^{i})$ & $
\left\{x,x+\frac{n}{2},3-i-x,3-i-x+\frac{n}{2}\right\}$ \\
  \hline
\end{tabular}\\
\end{center}
\noindent So we see that the equivalence class of $x\in[n]$  under the equivalence relation $C(\rho^{\frac{n}{2}})\vee
C(\tau\rho^{i})$ is
$$\overline{x}=
\left\{x,x+\frac{n}{2},3-i-x,3-i-x+\frac{n}{2}\right\},$$   so that $|\overline{x}|\in\{2,4\}.$ It is
not difficult to see that $|\overline{x}|=2$ if and only if either
$2x\equiv 3-i$ or $2x\equiv 3-i+\frac{n}{2}.$ Therefore
$b_{2}(C(\rho^{\frac{n}{2}})\vee C(\tau\rho^{i}))$ is $1$ if
$\frac{n}{2}$ is odd, $2$ if $\frac{n}{2}$ is even and $i$ is odd,
and $0$ if $\frac{n}{2}$ is even and $i$ is even.


\begin{thm}{\em Let $n \in \mathbb{N}_+$ be even.
According to whether $n$ is $0$ or $2$ mod $4$, the orbicycle
index polynomial of $D_{n}$ is given by
\begin{eqnarray*}
P_{D_{n}}^{\ orb}&=&P_{D_{n}}-P_{\mathbb{Z}_{n}}
+\frac{1}{2n}\left(
\sum_{d\mid\frac{n}{2}}\varphi
\left(\frac{n-1}{d},\frac{n}{2d}\right)x_{\frac{n}{d}}^{d}
+\left\{\begin{array}{cc} nx_{2}x_{4}^{\frac{n-2}{4}}\\
\frac{n}{2}x_{4}^{\frac{n}{4}}+\frac{n}{2}x_{2}^{2}x_{4}^{\frac{n-4}{4}}\end{array}\right.
\right)\\
&+&\frac{1}{n}\sum_{d\mid n}\varphi\left(\frac{n-2}{2d},\frac{n-1}{d},\frac{n}{d}\right)x_d^{\frac{n}{d}}
+\left(\frac{n+2}{2n}\right)\left(x_2^{\frac{n}{2}}+x_1^{2}x_2^{\frac{n-2}{2}}
+\left\{\begin{array}{cc}x_{2}x_{4}^{\frac{n-2}{4}}\\
x_{4}^{\frac{n}{4}}
\end{array}\right.
+\left\{\begin{array}{cc}
x_{2}x_{4}^{\frac{n-2}{4}}\\
x_{2}^{2}x_{4}^{\frac{n}{4}-1}
\end{array}\right.
\right).
\end{eqnarray*}
}\end{thm}
\begin{proof}
\begin{eqnarray*}
P_{D_{n}}^{\ orb}&=&\frac{1}{2n}\sum_{(\overline{g},h)\in
P(D_{n})}|\overline{g}|\prod_{k=1}^{n}x_{k}^{b_{k}(C(g)\vee
C(h))}\\
&=&P_{D_{n}}+\frac{1}{2n}\sum_{h\in
D_{n}}\prod_{k=1}^{n}x_{k}^{b_{k}(C(\rho^{\frac{n}{2}})\vee
C(h))}\\
&+&\frac{1}{n}\sum_{i=0}^{\frac{n-2}{2}}
\sum_{j=0}^{n-1}\prod_{k=1}^{n}x_{k}^{b_{k}(C({\rho}^{i})\vee
C({\rho}^{j}))}-P_{\mathbb{Z}_{n}}\\
&+&\frac{n+2}{4n}\sum_{h\in\{e,\tau,\rho^{\frac{n}{2}},
\rho^{\frac{n}{2}}\tau\}}\prod_{k=1}^{n}x_{k}^{b_{k}(C(\tau)\vee
C(h))}\\
&+&\frac{n+2}{4n}\sum_{h\in\{e,\rho\tau,\rho^{\frac{n}{2}},
\rho^{\frac{n}{2}+1}\tau\}}
\prod_{k=1}^{n}x_{k}^{b_{k}(C(\rho\tau)\vee C(h))}.
\end{eqnarray*}
We compute the last  four summands in the expression above

\begin{eqnarray*}
\sum_{h\in D_{n}}\prod_{k=1}^{n}x_{k}^{b_{k}(C(\rho^{\frac{n}{2}})\vee
C(h))}&=&
\sum_{i=0}^{n-1}\prod_{k=1}^{n}x_{k}^{b_{k}(C(\rho^{\frac{n}{2}})
\vee C(\rho^{i}))}
+\sum_{i=0}^{n-1}\prod_{k=1}^{n}x_{k}^{b_{k}
(C(\rho^{\frac{n}{2}}) \vee C(\tau\rho^{i}))}\\
&=&\sum_{i=0}^{n-1}x^{(\frac{n}{2},i)}_{\frac{n}{(\frac{n}{2},i)}}
+\sum_{i=0}^{n-1}\prod_{k=1}^{n}x_{k}^{b_{k}
(C(\rho^{\frac{n}{2}}) \vee C(\tau\rho^{i}))}\\
&=&\sum_{d\mid\frac{n}{2}}\varphi
\left(\frac{n-1}{d},\frac{n}{2d}\right)x_{\frac{n}{d}}^{d}
+\left\{\begin{array}{cc} nx_{2}x_{4}^{\frac{n-2}{4}}\\
\frac{n}{2}x_{4}^{\frac{n}{4}}+\frac{n}{2}x_{2}^{2}x_{4}^{\frac{n-4}{4}}
\end{array}\right.
\end{eqnarray*}

\begin{eqnarray*}
\sum_{i=0}^{\frac{n-2}{2}}\sum_{j=0}^{n-1}\prod_{k=1}^{n}
x_{k}^{b_{k}(C(\rho^{i})\vee
C(\rho^{j}))}&=&\sum_{i=0}^{\frac{n-2}{2}}\sum_{j=0}^{n-1}
x_{\frac{n}{(i,j,n)}}^{(i,j,n)}
=\sum_{d\mid n}\varphi\left(\frac{n-2}{2d},\frac{n-1}{d},\frac{n}{d}\right)x_d^{\frac{n}{d}}
\end{eqnarray*}

\begin{eqnarray*}
\sum_{h\in\{e,\tau,\rho^{\frac{n}{2}},
\rho^{\frac{n}{2}}\tau\}}\prod_{k=1}^{n}x_{k}^{b_{k}(C(\tau)\vee
C(h))}=2x_{2}^{\frac{n}{2}}+2\left\{\begin{array}{cc}x_{2}x_{4}^{\frac{n-2}{4}}\\
x_{4}^{\frac{n}{4}}
\end{array}\right.
\end{eqnarray*}

\begin{eqnarray*}
\sum_{h\in\{e,\rho\tau,\rho^{\frac{n}{2}},
\rho^{\frac{n}{2}+1}\tau\}}\prod_{k=1}^{n}x_{k}^{b_{k}(C(\rho\tau)\vee
C(h))}&=&2x_{1}^{2}x_{2}^{\frac{n-2}{2}}+2\left\{\begin{array}{cc}
x_{2}x_{4}^{\frac{n-2}{4}}\\
x_{2}^{2}x_{4}^{\frac{n}{4}-1}
\end{array}\right.\\
\end{eqnarray*}
\end{proof}

In this work we have extended P\'olya theory to the context of
orbiquotient sets. The main ingredient of the new theory is the
orbicycle index polynomial which we computed in various cases. We
expect that our construction will find applications in the study
of the topology of orbifolds and also in the theory of species. I
would be interesting to search for a further extension of P\'olya
theory within the context of rational combinatorics introduced in
\cite{BD1, BD2} based on the previous work \cite{RDEP} and further
discussed in \cite{BD3}. One should obtain a generalization of
P\'olya theory in which finite sets are replaced by finite
groupoids \cite{BaezDolan3}.

\subsection* {Acknowledgment} Thanks to Edmundo Castillo, Federico Hernandez, Eddy Pariguan, Sylvie Paycha, and
Domingo Quiroz.

\noindent hector.blandin@lacim.ca\\
\noindent Departamento de Matem\'aticas Puras y Aplicadas,\\
Universidad Sim\'on Bol\'ivar\\
Caracas, Venezuela\\

\noindent ragadiaz@gmail.com\\
\noindent Instituto de Matem\'{a}ticas y sus Aplicaciones\\
Universidad Sergio Arboleda\\
Bogot\'{a}, Colombia.

\end{document}